\pdfoutput=1

\RequirePackage[english]{babel}

\documentclass[a4paper,english]{amsart}

\usepackage[sc]{mathpazo}
\DeclareTextFontCommand{\textsl}{\fontfamily{ppl}\fontshape{sl}\selectfont}
\usepackage{eulervm}

\usepackage{microtype}
\usepackage{mathtools}
\usepackage{amssymb}
\usepackage{amsthm}
\usepackage[foot]{amsaddr}
\usepackage[all]{xy}
\usepackage{mathrsfs}
\usepackage{ifthen}
\usepackage{nicefrac}
\usepackage[
pdftitle={Splitting definably compact groups},
pdfauthor={Marcello Mamino},
pdfkeywords={Definable groups, o-minimality, fibre bundles},
colorlinks=false]{hyperref}

\makeatletter
\def\Hy@Warning#1{}
\makeatother

\makeatletter
\def\@setemails{% 
\ifnum\theg@author > 1 
\mbox{{\itshape E-mail addresses}:\space}{\ttfamily\emails}. 
\else 
\mbox{{\itshape E-mail address}:\space}{\ttfamily\emails}. 
\fi% 
}
\makeatother

\makeatletter
\def\ps@plain{\ps@empty
  \def\@oddfoot{\normalfont\normalsize \hfil\thepage\hfil}%
  \let\@evenfoot\@oddfoot}
\def\ps@firstpage{\ps@plain
  \def\@oddfoot{\normalfont\normalsize \hfil\thepage\hfil
     \global\topskip\normaltopskip}%
  \let\@evenfoot\@oddfoot
  \def\@oddhead{\@serieslogo\hss}%
  \let\@evenhead\@oddhead % in case an article starts on a left-hand page
}
\def\ps@headings{\ps@empty
  \def\@evenhead{%
    \setTrue{runhead}%
    \normalfont\normalsize
    \rlap{\thepage}\hfil
    \textsc{\lsstyle\MakeLowercase{\shortauthors}\hfil}}%
  \def\@oddhead{%
    \setTrue{runhead}%
    \normalfont\normalsize \hfil
    \textsc{\lsstyle\MakeLowercase{\rightmark{}{}}}\hfil\llap{\thepage}}%
  \let\@mkboth\markboth
}
\makeatother
\makeatletter
\let\oldupchars\upchars@
\def\upchars@{\oldupchars\def\-{\U-}}
\makeatother

\def\mathshift{$}
\catcode`$=13
\def${\ifinner\expandafter\mathshift\else\expandafter\myshift\fi}
\def\myshift#1${\raisebox{0ex}[0ex][0ex]{\mathshift#1\mathshift}}
\AtBeginDocument{\catcode`$=13\iffalse$\fi}

\makeatletter
\let\oldsection\section
\def\newsection#1{\oldsection{\lsstyle #1}}
\def\newsectionr#1{\oldsection*{\lsstyle #1}}
\def\section{\@ifstar\newsectionr\newsection}
\makeatother

\makeatletter
\newtheorem{ghost@theorem}{}[section]
\def\@maketheorem#1=#2;{
	\newtheorem{#1}[ghost@theorem]{#2}}
\def\maketheorem#1{
	\@for\@x:=#1\do{
		\expandafter\@maketheorem\@x;}}
\makeatother

\makeatletter
\def\os@arabic#1{\oldstylenums{\number #1}}
\def\osarabic#1{\texorpdfstring
			{\expandafter\os@arabic\csname c@#1\endcsname}
			{\expandafter\@arabic\csname c@#1\endcsname}}

\@ifundefined{thepart}{}{}
\@ifundefined{thechapter}{%
	\def\thesection    {\osarabic{section}}

}{%
	\def\thechapter    {\osarabic{chapter}}

	\def\thesection    {\thechapter.\osarabic{section}}
}

\def\theghost@theorem{\thesection.\osarabic{ghost@theorem}}

\def\thepage{\osarabic{page}}

\makeatother

\def\myscriptsize{\tiny}

\def\rest#1{\mkern-5mu\upharpoonright\mkern-5mu\mathchoice
	{\raisebox{-.5ex}{$\mathsurround=0pt\scriptstyle{#1}$}}
	{\raisebox{-.5ex}{$\mathsurround=0pt\scriptstyle{#1}$}}
	{\raisebox{-.3ex}{$\mathsurround=0pt\scriptscriptstyle{#1}$}}
	{\raisebox{-.3ex}{$\mathsurround=0pt\scriptscriptstyle{#1}$}}}

\def\bundle#1{\mathscr #1}
\let\covering\calligraphic
\let\family\calligraphic
\def\onthebase#1{\overline{#1}}
\def\indmap#1{\check{#1}}
\def\eqdef{\stackrel{\textsl{\tiny def}}{=}}
\def\I{\mathrm I}
\def\T{\mathrm T}
\def\N{\mathbb N}
\def\R{\mathbb R}
\def\M{\mathcal M}
\def\Mst{\mathrm R}
\def\i{\mathbf i}
\def\j{\mathbf j}
\def\k{\mathbf k}
\def\SU{\operatorname{SU}}
\def\SO{\operatorname{SO}}
\def\Ralg{\R^\textsl{\myscriptsize alg}}
\def\wildcard{\text{--}}
\def\st{\textrm{ : }}
\def\closure#1{\overline{#1}}
\def\boundary#1{\delta #1}
\def\dom#1{\operatorname{Dom}(#1)}

\def\id{\operatorname{Id}}
\def\P{$\mathscr P$}
\def\identity#1{\mathrm{e}_{\scriptscriptstyle{#1}}}
\def\pidef{\pi^{\textsl{\myscriptsize def}}}
\def\gcenter#1{\operatorname{Z}(#1)}
\def\cccenter#1{\operatorname{Z}^0(#1)}
\def\derived#1{[\mkern1mu#1,#1\mkern1mu]}
\def\eqvdots{\mathmakebox[\widthof{${}={}$}][c]{\vdots}}
\def\quotient#1#2{\mathchoice
	{#1/\raisebox{-.5ex}{$\mathsurround=0pt\displaystyle #2$}}
	{#1/\raisebox{-.5ex}{$\mathsurround=0pt\textstyle #2$}}
	{#1/\raisebox{-.3ex}{$\mathsurround=0pt\scriptstyle #2$}}
	{#1/\raisebox{-.1ex}{$\mathsurround=0pt\scriptscriptstyle #2$}}}
\def\intervalcc#1#2{[\mkern1mu#1,#2\mkern1mu]}
\def\intervalco#1#2{[\mkern1mu#1,#2\mkern-1mu)}
\def\intervaloc#1#2{(\mkern-1mu#1,#2\mkern1mu]}
\def\intervaloo#1#2{(\mkern-1mu#1,#2\mkern-1mu)}

\let\n\oldstylenums
\def\mydate#1#2#3{\hbox{\n{#1}$\cdot${\sc#2}$\cdot$\n{#3}}}
\def\new#1{\textsc{#1}}
\def\-{\nobreakdash-\hspace{0pt}}
\def\U-{\raise0.2ex\hbox{-}}
\def\url#1{\href{#1}{url\nobreakdash---\texttt{#1}}}
\def\doi#1{\href{http://dx.doi.org/#1}{\texttt{doi:#1}}}
\def\arXiv#1{\href{http://arxiv.org/abs/#1}{\texttt{arXiv:#1}}}
\def\mailto#1{\href{mailto:#1}{\texttt{#1}}}

\theoremstyle{plain}
\maketheorem{%
	theorem=Theorem,%
	lemma=Lemma,%
	corollary=Corollary,%
	fact=Fact%
}
\theoremstyle{definition}
\maketheorem{%
	definition=Definition,%
	example=Example%
}
\theoremstyle{remark}
\maketheorem{%
	observation=Observation,%
	remark=Remark%
}

\begin{document}

\title[Splitting definably compact groups]
	{Splitting definably compact groups\\
	in o\-minimal structures}
\author[M. Mamino]{\lsstyle Marcello Mamino}
\address{Classe di Scienze -- Scuola Normale Superiore,
	Piazza dei Cavalieri \n{7},
	\n{56126} Pisa,
	Italy.}
\email{\mailto{m.mamino@sns.it}}
\date{\mydate{23}{xii}{2009}, revised: \mydate{22}{ix}{2010}}

\def\abstractname{\lsstyle Abstract}
\begin{abstract}
An argument of A.~Borel~\cite[Proposition~\n{3}.\n{1}]{Bo61} shows that
every compact connected Lie group is homeomorphic to the Cartesian product
of its derived subgroup and a torus. We prove a parallel result for
definably compact definably connected groups definable in an o\-minimal
expansion of a real closed field. As opposed to the Lie case, however, we
provide an example showing that the derived subgroup may not 
have a definable semidirect complement.
%be a semidirect
%factor of the group.
\end{abstract}
\maketitle

\section{Introduction}\label{sect-intro}

\noindent It is known since~\cite{Pi88} that groups definable in o\-minimal structures
are \emph{ipso facto} topological groups and definable manifolds, moreover
in a unique way. The (definable) homotopic invariants of such groups, in
particular in the definably compact case, have been the subject
of several recent papers. Using the homotopic information
gained so far, we are able to
show an analogue of~\cite[Proposition~\n{3}.\n{1}]{Bo61} for
definably compact definably connected groups definable in
an o\-minimal structure expanding a real closed field.
Namely, every such group is definably
homeomorphic to the Cartesian product
%characterize
%(theorem~\ref{th-groups-are-products})
%the topology of a definably
%connected definably
%compact group as the topological product
of its derived subgroup -- which is definable by~\cite[Corollary~\n{6}.\n{4}]{HruPePi08B} -- and a
definable Abelian group---which, precisely, is the quotient of the
connected component of the identity of the group's center modulo its -- finite --
intersection with the derived subgroup
(Theorem~\ref{th-groups-are-products}).

For compact Lie groups a stronger result is known---as the
Borel\-Scheerer\-Hofmann splitting theorem:
see~\cite[Theorem~\n{9}.\n{39}]{Hofmann-Morris}. The derived subgroup of a
compact Lie group is indeed a semidirect factor of the group.
In Example~\ref{ex-no-semidirect-product}, we show that this may not be
the case in the definable context.

Proposition~\n{3}.\n{1} in~\cite{Bo61} -- see
also~\cite[Fact~\n{4}.\n{3}]{Ba09} -- is the Lie counterpart of our
main result (Theorem~\ref{th-groups-are-products}). It is proven by
induction on the dimension: roughly
\begin{quotation}
Let~$G$ be our Lie group and let~$\cccenter{G}$ denote
the connected component of the identity of the center of~$G$. Take a
connected subgroup~$Z$ of~$\cccenter{G}$ having codimension~$1$.
By induction,
we can assume the statement on~$Z\cdot\derived{G}$,
where $\derived{G}$ denotes the derived subgroup.
The base space of the principal fibre bundle induced by the quotient
$\quotient{G}{Z\cdot\derived{G}}$
is homeomorphic to a
circle, hence, using~\cite[Corollary~\n{18}.\n{6}]{Steenrod}, the bundle is trivial.
\end{quotation}
However simple,
this technique can not be adapted blindly to the o\-minimal
setting. The problem resides in the inductive step, which is possible
because
any compact Abelian Lie group factorizes in a product of
$1$\-dimensional groups. The same is known not to hold for definably
compact groups --
see~\cite[Section~\n{5}]{PeSt99} -- at least if by~\emph{factorizes in a
product of}  we
mean~\emph{is definably isomorphic to the direct product of,} while it is unknown if we
mean~\emph{is definably homeomorphic to the Cartesian product of.}

Here we present a different approach, based on the study of the homotopic
properties of definable fibre bundles. The notion of definable fibre
bundle has been introduced in~\cite{BeMaOt08}, in order to characterize the
higher homotopy groups of definably compact groups. There, it is shown that
quotients of definable groups induce definable fibre
bundles (Lemma~\n{2}.\n{2}) and that
definable fibre bundles have the homotopy lifting property
(Theorem~\n{2}.\n{3}), hence they share many of the
properties that well behaved fibrations are supposed to enjoy---for example,
by~\cite[Theorem~\n{4}.\n{9}]{BaOt09}, we have the usual
long exact sequence of definable homotopy groups.
Further on this road, we give a definition
of~\emph{definable bundle map}. Then we show that, under suitable
assumptions,
our bundle maps behave as their topological counterparts.
As a consequence, we can give a sufficient condition for a definable fibre
bundle, with definably connected fibre, on a torus -- where \emph{torus}
means $\quotient{{\intervalcc{0}{1}}^n}{\sim}$ -- to have a definable
section.
% As a consequence, we can give a
% necessary and sufficient condition for a \emph{principal} definable fibre
% bundle on a torus -- where \emph{torus} means $\quotient{{\intervalcc{0}{1}}^n}{\sim}$ -- to be trivial.
Our condition
is purely homotopic. (Elaborating on the same lines, it is possible to
give a necessary and sufficient condition for a \emph{principal} definable
fibre bundle on a torus to be trivial. This has been done in the author's
Ph.D.\ thesis.)
%---in fact, it \emph{would} be equivalent to the
%vanishing of
%all the relevant obstructions, if we \emph{had} an obstruction
%theory in the o\-minimal context (which should not be hard to construct).
Going back to Borel's argument,
knowing that any definably compact definably connected Abelian group is
definably homotopy equivalent to a torus -- Theorem~\n{3}.\n{4} in~\cite{BeMaOt08}
-- we are able to substitute a factorized torus for the potentially wild
center of~$G$, thus solving the difficulty.

\subsection*{Conventions and notations}

Fix~$\M\eqdef(\Mst,<,0,1,+,\cdot,\dotsc)$, an o\-minimal expansion of a real closed field. Throughout this
paper~\new{definable} will mean~\emph{definable in~$\M$,} with parameters.
We put on~$\Mst$
the order topology and on~$\Mst^n$ the product topology. On the subsets of
them, we put the subset topology. We will assume some knowledge of
o\-minimality and, in particular, of groups definable in o\-minimal
structures---for references see~\cite{VanDenDries} and~\cite{Ot08}
respectively.

Additionally, for us, \new{definable map} means~\emph{continuous definable function},
and~\new{definable path} means~\emph{continuous definable curve}, i.e.\ function having
$\I\eqdef\intervalcc{0}{1}\subset\Mst$ as domain. Clearly \emph{functions} and~\emph{curves} are not
necessarily continuous. A~\new{definable homotopy} is a definable map~$f$
having~$\begin{smallmatrix}\textup{\tiny something}\\
\textup{\tiny definable}\end{smallmatrix}\mkern-3mu\times\I$ as domain,
in this situation $f_t(\wildcard)$ is short
for~$f(\wildcard,t)$; definable homotopic properties, such as~\emph{definable
homotopy equivalence,} are specified in terms of definable homotopies as in
the topological setting---for more details see~\cite{BaOt09}.

\subsection*{Acknowledgements}

\hyphenation{Al-es-san-dro Ber-ar-duc-ci An-na-li-sa Con-ver-sa-no}
The author would like to express his gratitude to Alessandro Berarducci
and Annalisa Conversano for the helpful conversations had during the
preparation of the present paper.

\section{Definable fibre bundles}\label{sect-on-bundles}

\noindent In this section we concentrate on~\emph{definable bundle maps}
and~\emph{homotopies} of them.
Most of the work is an adaptation of the usual topological arguments to
the definable context. In particular, virtually all the o\-minimality
here is embedded in Lemma~\ref{th-fam-inverse} and
Observation~\ref{th-slicing-lemma}. The decision not to provide our
bundles with a structure group has been a considerate one, since both
alternatives would have worked, for our purpose. Basically we decided for
the shortest route, however this often forced us to assume the fibre of
our bundles to be~\emph{locally definably compact.} The necessity of such
assumption is shown by Example~\ref{ex-locally-compact}, its reason resides
in Lemma~\ref{th-fam-inverse}.
%---which, rather informally, says that the group of all the
%homeomorphisms of a locally definably compact definable set, with the
%compact open topology, is a topological group.

\begin{definition}\label{def-bundle}
A \new{definable fibre bundle}~$\bundle{B}=(B,X,p,F)$ consists
of a definable set~$F$ and a definable map~$p\colon B\to X$ between definable sets
$B$ and~$X$ with the following properties:
\begin{itemize}
\item[(i)] $p$ is onto
\item[(ii)] there is a finite covering by open definable subsets~$\covering{U}=\{U_i\}_i$ of~$X$ such that for each~$i$
there is a definable homeomorphism~$\phi_i\colon U_i\times F\to p^{-1}(U_i)$.
\end{itemize}
The set~$B$ is said to be the~\new{bundle space} of~$\bundle{B}$, the
set~$X$ is said to be its~\new{base}, and~$F$ its~\new{fibre}. Any
covering as in~(ii) is said to be a~\new{trivialization covering} and the
maps~$\phi_i$ are called~\new{trivialization maps}. We will refer to the
map~$p$ as the~\new{projection}, moreover, fixed a trivialization
covering~$\covering{U}=\{U_i\}_i$ with the associated
trivialization maps, we will call~\new{projection on the fibre} each of the
maps~$p_i\colon p^{-1}(U_i)\to F$ sending $x\in p^{-1}(U_i)$ to the second component
of~$\phi_i^{-1}(x)$.
Finally we will call~\new{definable cross section} of~$\bundle{B}$ any definable
map~$s\colon X\to B$ such that $p\circ s$ is the identity on~$X$.
\end{definition}

\begin{definition}\label{def-bmap}
Let $\bundle{B}=(B,X,p,F)$
and~$\bundle{B'}=(B',X',p',F)$ be definable fibre bundles
having the same fibre. A
\new{definable bundle map}~$f\colon\bundle{B}\to\bundle{B'}$ is a
definable map~$f\colon B\to B'$ such that:
\begin{itemize}
\item[(i)] there is a definable map~$\onthebase{f}\colon X\to X'$ such
that the following diagram commutes
\[\xymatrix{
B \ar[r]^f \ar[d]^p & B\mathrlap{\smash{'}} \ar[d]^{p\mathrlap{\smash{'}}}\\
X \ar[r]^{\onthebase{f}} & X\mathrlap{\smash{'}}
}\]
\item[(ii)] for each~$x\in X$ the map~$f\rest{p^{-1}(x)}$ is an
homeomorphism onto~$p'^{-1}\circ\onthebase{f}(x)$.
\end{itemize}
The map~$\onthebase{f}$ is said to be induced by $f$ on the base spaces.
\end{definition}

\begin{definition}\label{def-bundle-equiv}
A definable bundle map~$f\colon\bundle{B}\to\bundle{B'}$ between definable fibre
bundles $\bundle{B}$ and~$\bundle{B'}$ is said to be a \new{definable bundle
isomorphism} whenever it is an homeomorphism (of the bundle spaces of
$\bundle{B}$ and~$\bundle{B'}$). In this case $\bundle{B}$
and~$\bundle{B'}$ are said to be~\new{isomorphic}.
\end{definition}

Some examples of definable fibre bundles are Cartesian products of definable
sets, definable covers, and quotients of definable groups (for the latter
see Fact~\ref{th-quotient-bundle} below).
As usual, we will call~\new{trivial} a definable fibre bundle which is
isomorphic to a Cartesian product. The following lemma gives us a useful
criterion to tell when a definable bundle map is an isomorphism, under the
fairly general hypothesis that the fibre is~\new{locally definably
compact}---i.e.\ that each point of it has a definably compact definable neighbourhood.

\begin{lemma}\label{th-equivalence}
Let $f\colon\bundle{B}\to\bundle{B'}$ be a definable bundle map between
the definable fibre bundles $\bundle{B}=(B,X,p,F)$
and~$\bundle{B'}=(B',X',p',F)$, with $F$ locally definably
compact.  Suppose that the map $\onthebase{f}\colon X\to X'$ induced by~$f$ on the base
spaces is an homeomorphism.  Then $f$ is a definable bundle isomorphism.
\end{lemma}

We will postpone the proof of Lemma~\ref{th-equivalence} until some
additional machinery has been developed.

\begin{definition}\label{def-family}
Let $X$, $Y$, and~$Z$ be definable sets. We call a
family~$\family{F}=\{f_x\}_{x\in X}$ of functions from $Y$ to
$Z$ a \new{definable family} if the function
\begin{align*}
f\colon X\times Y &\to Z \\
(x,y)&\mapsto f_x(y)
\end{align*}
%f\colon X\times Y\ni(x,y)\mapsto f_x(y)\in Z\]
is definable. We call $\family{F}$ a
\new{continuous family} whenever $f$ is continuous.
\end{definition}

\begin{lemma}\label{th-fam-inverse}
Let $X$ and~$Y$ be definable sets, with $Y$ locally definably compact. Consider a
definable continuous family~$\family{F}=\{f_x\}_{x\in X}$
\emph{of homeomorphisms} of~$Y$. We claim that the definable
family~$\family{F}^{-1}\eqdef\{f_x^{-1}\}_{x\in X}$ of
homeomorphisms of~$Y$, which we will call the~\new{inverse}
of~$\family{F}$, is itself continuous.
\end{lemma}
\begin{proof}
Clearly $\family{F}^{-1}$ is definable: we want to prove that the function
\begin{align*}
g\colon X\times Y &\to Y \\
(x,y)&\mapsto f_x^{-1}(y)
\end{align*}
%g\colon X\times Y\ni(x,y)\mapsto f_x^{-1}(y)\in Y\]
is continuous. By contradiction, suppose $g$ to be discontinuous
at~$(a,b')\in X\times Y$ with~$b'=f_a(b)$ for some~$b\in Y$.
Hence, there are definable paths $\gamma_1\colon\I\to X$
and~$\gamma_2\colon\I\to Y$ such that $\gamma_1(0)=a$ and~$\gamma_2(0)=b'$
but the curve~$x\mapsto f^{-1}_{\gamma_1(x)}\circ\gamma_2(x)$ does not
converge to~$b$ for $x$ going to~$0$.
Now, for each~$x\in\I$ define the curve
\begin{align*}
\gamma'_x\colon\I &\to Y \\
t &\mapsto
%\gamma'_x\colon\I\ni t\mapsto
	\begin{cases}
		\gamma_2\left((1-2t)x\right) &\text{for~$t<\frac{1}{2}$}\\
		f_{\gamma_1\left((2t-1)x\right)}(b) &\text{for~$\frac{1}{2}\le t$}
	\end{cases}
\end{align*}
which is continuous, definable, and joins $\gamma'_x(0)=\gamma_2(x)$
to~$\gamma'_x(1)=f_{\gamma_1(x)}(b)$.
Fix a definably compact definable neighbourhood~$\closure{V}\subset Y$ of~$b$,
and let~$t\colon\I\to\I$ be the function mapping $x$ to the least~$t$ such
that~$f^{-1}_{\gamma_1(x)}\circ\gamma'_x(t)$ is in~$\closure{V}$, which is
definable, and is well defined observing that
\[f^{-1}_{\gamma_1(x)}\circ\gamma'_x(1)=b\in\closure{V}\]
Define the curve, not necessarily continuous
\begin{align*}
\gamma_3\colon\I &\to \closure{V}\\
x &\mapsto f^{-1}_{\gamma_1(x)}\circ\gamma'_x\circ t(x)
\end{align*}
%\gamma_3\colon\I\ni x\mapsto f^{-1}_{\gamma_1(x)}\circ\gamma'_x\circ
%t(x)\in\closure{V}
Now, the limit~$b''\eqdef\lim_{x\to0}\gamma_3(x)$ exists because
$\closure{V}$ is definably compact.
Also, for each~$x\in\I$, either $t(x)=0$ or $\gamma_3(x)$ lies on the
boundary of~$\closure{V}$ (inside~$Y$), hence, by
o\-minimality,
there is an~$\epsilon>0$ such that either of these conditions holds for
all~$x\in\intervaloo{0}{\epsilon}$. In both cases~$b''\neq b$: if $t(x)=0$
for all~$x\in\intervaloo{0}{\epsilon}$, then, for
$x\in\intervaloo{0}{\epsilon}$, we have that $\gamma_3(x) =
f^{-1}_{\gamma_1(x)}\circ\gamma_2(x)$, which does not converge to~$b$ by
hypothesis. Otherwise $b''$ must lie
on the boundary of~$\closure{V}$, hence can not equal~$b$.
%Clearly, for each~$x\in\I$, either $t(x)=0$ or $\gamma_3(x)$ lies on the
%boundary of~$\closure{V}$, hence exists~$b''\eqdef\lim_{x\to0}\gamma_3(x)$
%and~$b''\neq b$.
However
\[
f_a(b'')	=\lim_{x\to0}f_{\gamma_1(x)}\circ\gamma_3(x)
		=\lim_{x\to0}\gamma'_x\circ t(x)
		=b'=f_a(b)
\]
where the second last equality holds observing that, for $x$ going
to~$0$, the paths~$\gamma'_x$ converge uniformly to the constant path at~$b'$.
This contradicts $f_a\colon Y\to Y$ being an homeomorphism, and, in
particular, \n{1}~to~\n{1}.
\end{proof}

The following example shows that in the statement of
Lemma~\ref{th-fam-inverse} we must assume $Y$~to be locally definably
compact. Similarly, we must assume the fibre to be locally definably
compact in Lemma~\ref{th-equivalence}.

\begin{example}\label{ex-locally-compact}
Working in the o\-minimal structure~$\Ralg$, let $Q$ be the definable
set~$\R^{>0}\times\R^{>0}\cup\{(0,0)\}\subset\R^2$.
We will give a continuous definable family~$\family{F}$ of homeomorphisms
of~$Q$ whose inverse~$\family{F}^{-1}$ is not continuous.

For $t\in\R^{\ge0}$, let $\alpha_t\colon\R^{>0}\to\R^{>0}$ be the map
\[
\alpha_t(x) = \begin{cases}
	\min\left(1,t+\left(\frac{x-t}{t}\right)^2\right)
		&\text{if $t>0$}\\
	1	&\text{if $t=0$}
\end{cases}
\]
and define the family~$\family{F}=\{f_t\}_{t\in\R^{\ge0}}$ as
\[
f_t(x,y) = \begin{cases}
	(0,0)	&\text{if $(x,y)=(0,0)$}\\
	\alpha_t\left(\frac{\raisebox{2pt}{$\mathsurround=0pt\scriptstyle{y}$}}{x}\right)\cdot(x,y)
		&\text{otherwise}
\end{cases}
\]
which is clearly continuous (at~$(0,0)$ because
$\alpha_\wildcard(\wildcard)$ is bounded by $1$). However, $\family{F}^{-1}$ is
not continuous, in fact, considering the path~$\gamma\colon\I\ni
t\mapsto(t,t^2)\in Q$, we have
\[
f_0^{-1}\circ\gamma(0)=(0,0)\neq(1,0)=\lim_{t\to0}f_t^{-1}\circ\gamma(t)
\]
Using $\family{F}$, an example can be constructed of definable bundle map
which
is not a definable bundle isomorphism,
but induces an homeomorphism on the base spaces.
In fact, let $\bundle{B}$ be the trivial bundle~$\R^{\ge0}\times Q$ with
base~$\R^{\ge0}$ and fibre~$Q$. The map
\begin{align*}
	f\colon\bundle{B}&\to\bundle{B} \\
		\left(t,\left(x,y\right)\right) &\mapsto
		\left(t,f_t\left(x,y\right)\right)
\end{align*}
is a definable bundle map and enjoys the claimed property.
\end{example}

\begin{proof}[Proof of Lemma~\ref{th-equivalence}]
Clearly $f$ is \n{1}~to~\n{1}, so suffices to show that $f^{-1}$ is
continuous. By contradiction, let $f(x)$ be a
discontinuity point of~$f^{-1}$. Fix definable trivialization coverings
$\covering{U}=\{U_i\}_i$
and~$\covering{U'}=\{U'_j\}_j$ for~$\bundle{B}$
and~$\bundle{B'}$ respectively, and fix $i$ and~$j$ such that $p(x)\in U_i$
and $p'\circ f(x)\in U'_j$. Observe that the set $U\eqdef U_i\cap\onthebase{f}^{-1}(U'_j)$
is an open neighbourhood of~$p(x)$, so $\onthebase{f}(U)$ is an open
neighbourhood of~$p'\circ f(x)=\onthebase{f}\circ p(x)$, and $p'^{-1}\circ\onthebase{f}(U)$ is an
open neighbourhood of~$f(x)$. Hence suffices to show
that $f^{-1}\rest{p'^{-1}\circ\onthebase{f}(U)}$ is continuous. For, let
$\phi_i$ and~$\phi'_j$ be trivialization maps associated to~$U_i$
and~$U'_j$ respectively. Since $\phi'_j$ sends
$\onthebase{f}(U)\times F$ isomorphically
onto~$p'^{-1}\circ\onthebase{f}(U)$, we can reduce to prove that
$\phi_i^{-1}\circ f^{-1}\circ\phi_j\rest{\onthebase{f}(U)\times F}$ is
continuous. However, this is equivalent to the
family
\[\family{F}^{-1}\eqdef\left\{p_i\circ
f^{-1}\circ\phi_j(x,\wildcard)\right\}_{x\in\onthebase{f}(U)}
\]
of
homeomorphisms of the fibre being continuous, which we have by
Lemma~\ref{th-fam-inverse}, since it is the inverse
of
\[\family{F}\eqdef\left\{p'_j\circ
f\circ\phi_i\left(\onthebase{f}^{-1}(x),\wildcard\right)\right\}_{x\in\onthebase{f}(U)}
\]
where $p_i$ and~$p'_j$ are the projections on the fibre corresponding
to~$\phi_i$ and~$\phi'_j$.
\end{proof}

\begin{definition}\label{def-ind-bundle}
Let $\bundle{B}=(B,X,p,F)$ be a definable fibre bundle, and
let $f\colon Y\to X$ be a definable map from the definable set~$Y$ to the
base space of~$\bundle{B}$. We define the~\new{induced
bundle}~$f^{-1}(\bundle{B})\eqdef(A,Y,p',F)$ where
\begin{itemize}
\item[(i)] $A=\{(x,y)\in B\times Y \st p(x)=f(y)\}$
\item[(ii)] $p'\colon A\ni(x,y)\mapsto y\in Y$
\end{itemize}
Moreover we will refer to the map~$\indmap{f}\colon A\ni(x,y)\mapsto x\in
B$ as induced by~$f$.
\end{definition}

\begin{lemma}\label{th-ind-bundle}
Using the same notations of Definition~\ref{def-ind-bundle} we have that
$f^{-1}(\bundle{B})$ is actually a definable fibre bundle, and
$\indmap{f}$ is a definable bundle map~$f^{-1}(\bundle{B})\to\bundle{B}$.
\end{lemma}
\begin{proof}
Straightforward: $p'$ is clearly onto, so suffices to fix a trivialization
covering~$\covering{U}=\{U_i\}_i$ for~$\bundle{B}$ and check
that~$\covering{V}\eqdef\{f^{-1}(U_i)\}_i$ is a trivialization
covering of~$f^{-1}(\bundle{B})$. For each~$i$, let $\phi_i$ be a trivialization
map for~$\bundle{B}$ associated to~$U_i$, and observe that the following
\begin{align*}
\phi'_i\colon f^{-1}(U_i)\times F &\to p'^{-1}\circ f^{-1}(U_i) \\
(x,y) &\mapsto \left(\phi_i\left(f(x),y\right),x\right)
\end{align*}
%\ni(x,y)\mapsto\left(\phi_i\left(f(x),y\right),x\right)\in
%p'^{-1}\circ f^{-1}(U_i)
is an homeomorphism, hence we may take it as the trivialization map for~$f^{-1}(\bundle{B})$
associated to~$f^{-1}(U_i)$. That~$\indmap{f}$ is a bundle map follows
immediately by inspection of the definition.
\end{proof}

\begin{observation}\label{th-composition-ind}
Consider three definable sets~$X_1$,$X_2$, and~$X_3$ with two definable
maps~$f_1\colon X_1\to X_2$ and~$f_2\colon X_2\to X_3$ between them. Take a definable
fibre bundle~$\bundle{B}=(B,X_3,p,F)$, and
define~$\bundle{B'}=f_1^{-1}(f_2^{-1}(\bundle{B}))$
and~$\bundle{B''}=g^{-1}(\bundle{B})$, where $g=f_2\circ f_1$.
It is easy to see from Definition~\ref{def-ind-bundle} -- and also,
when~$F$ is locally definably compact, follows
from Theorem~\ref{th-equivalence-theo} below --
that $\bundle{B'}$ and~$\bundle{B''}$ are isomorphic. In addition, there
is a definable bundle isomorphism~$h\colon\bundle{B'}\to\bundle{B''}$ such that
$\indmap{g}\circ h=\indmap{f}_2\circ\indmap{f}_1$.
\end{observation}

\begin{theorem}\label{th-equivalence-theo}
Let $\bundle{B}=(B,X,p,F)$
and~$\bundle{B'}=(B',X',p',F)$ be definable fibre bundles
having the same fibre~$F$, and let~$f\colon\bundle{B}\to\bundle{B'}$ be a
definable bundle map between them.
Suppose that $F$ is a locally definably compact definable set.
Then the the bundle~$\bundle{B''}\eqdef\onthebase{f}^{-1}(\bundle{B'})$,
where~$\onthebase{f}$ is the map induced by~$f$ on the base spaces, is
isomorphic to~$\bundle{B}$. Moreover there is a definable bundle
isomorphism~$g\colon\bundle{B}\to\bundle{B''}$ such
that~$\indmap{\onthebase{f}}\circ g=f$.
\end{theorem}
\begin{proof}
By definition, $\bundle{B''}=(B'',X,p'',F)$ where
\begin{align*}
	B''&=\left\{(x,y)\in B'\times X \st
	p'(x)=\onthebase{f}(y)\right\}\\
	p''&\colon(x,y)\mapsto y
\end{align*}
We claim that
\begin{align*}
g\colon B &\to B'' \\
x &\mapsto\left(f(x),p(x)\right)
\end{align*}
works. In fact, $g$ is clearly a definable bundle map, moreover it induces
the identity on the common base space~$X$ of~$\bundle{B}$
and~$\bundle{B''}$, hence, by Lemma~\ref{th-equivalence}, it is a
definable bundle isomorphism. The identity~$\indmap{\onthebase{f}}\circ
g=f$ is immediate.
\end{proof}

Observe that the hypothesis of local compactness can not be removed from
Theorem~\ref{th-equivalence-theo}, in fact any definable bundle map
for which the thesis of Lemma~\ref{th-equivalence} fails would make
Theorem~\ref{th-equivalence-theo} fail as well.

\begin{definition}
Let $\bundle{B}=(B,X,p,F)$
and~$\bundle{B'}=(B',X',p',F)$ be definable fibre bundles
having the same fibre, and let $f$ and~$g$ be definable bundle
maps~$\bundle{B}\to\bundle{B'}$. We call~\new{homotopy of definable bundle
maps} between $f$ and~$g$ a definable map~$h\colon B\times\I\to B'$
having the following properties:
\begin{itemize}
\item[(i)] $h_0\eqdef h(\wildcard,0)=f$
\item[(ii)] $h_1\eqdef h(\wildcard,1)=g$
\item[(iii)] For each $t\in\I$ the map~$h_t\eqdef h(\wildcard,t)$ is a
definable bundle map.
\end{itemize}
In this situation, the function~$\onthebase{h}\colon
X\times\I\ni(x,t)\mapsto\onthebase{h_t}(x)\in X'$ is a definable homotopy
between~$\onthebase{f}$ and~$\onthebase{g}$, to which we will refer as
the~\new{homotopy induced on the base spaces} by~$h$.

Moreover, we will say that $h$ is~\new{stationary with the induced
homotopy}~$\onthebase{h}$ if for any subinterval~$\intervalcc{t_1}{t_2}$ of~$\I$ and
any~$x\in B$ the following happens: $h\rest{\{x\}\times\intervalcc{t_1}{t_2}}$ is
constant if and only if
$\onthebase{h}\rest{\{p(x)\}\times\intervalcc{t_1}{t_2}}$ is constant.
\end{definition}

The following fact is Lemma~\n{2}.\n{10} in \cite{BeMaOt08}.

\begin{fact}\label{th-slicing-lemma}
Let $X$ and~$B$ be definable sets, and consider a definable homotopy~$f\colon X\times\I\to B$.
Let $\covering{U}$ denote a finite covering by open definable subsets of~$B$. Then, there
are finitely many definable
maps~$0\equiv g_1\le\dotsb\le g_{k+1}\equiv1\colon X\to\I$
such that, for each~$x\in X$ and
each~$i\in\{1,\dotsc,k\}$, the set~$f(x,\intervalcc{g_i(x)}{g_{i+1}(x)})$
is entirely contained in some element of~$\covering{U}$.
\end{fact}

\begin{theorem}\label{th-covering-homotopy}
Let $\bundle{B}=(B,X,p,F)$
and~$\bundle{B'}=(B',X',p',F)$ be definable fibre bundles
having the same fibre. Let $f\colon\bundle{B}\to\bundle{B'}$ be a
definable bundle map and $\tilde{h}\colon X\times\I\to X'$ be a definable homotopy
of~$\onthebase{f}$, i.e.\ $\tilde{h}_0=\onthebase{f}$. Then
$\tilde{h}$ can be lifted to an homotopy~$h$ of definable bundle maps such
that $h_0=f$ and $\onthebase{h}=\tilde{h}$. Moreover, $h$ can be chosen so
that it is stationary with~$\onthebase{h}$.
\end{theorem}
\begin{proof}
Let~$\covering{U}=\{U_i\}_{1\le i\le n}$ be a trivialization
covering for~$\bundle{B'}$, and take a definable shrinking~$\covering{U'}=\{U'_i\}_{1\le i\le n}$
of~$\covering{U}$; i.e.\ $\covering{U'}$ covers~$X'$ and for
each~$n$ holds~$\closure{U'_n}\subset U_n$ (it exists
by~\cite[Chapter~\n{6} Lemma~\n{3}.\n{6}]{VanDenDries}).
By Fact~\ref{th-slicing-lemma} we have
finitely many definable maps~$0\equiv g_1\le\dotsb\le g_{k+1}\equiv
1\colon X\to\I$ such that for each~$x\in X$ and for
each~$j\in\{1,\dotsc,k\}$ the
set~$\tilde{h}(x,\intervalcc{g_j(x)}{g_{j+1}(x)})$ is entirely
contained in some element of~$\covering{U'}$. 
%Not to clutter the argument with quantifications, from now on, we will assume an
%implicit \emph{for each~$i\in\{1,\dotsc,n\}$ and each~$j\in\{1,\dotsc,k\}$},
%unless stated otherwise.

Let
\[
V_{j,i}=\left\{x\st\tilde{h}\left(x,\intervalcc{g_j(x)}{g_{j+1}(x)}\right)\subset
U'_i\right\}
\]
Clearly, for each~$j$ the finite
family~$\{V_{j,i}\}_i$ of open definable sets is
a covering of~$X$. Choose definable open subsets
$W_{j,i}$ of~$X$ so that~$\closure{W_{j,i}}\subset V_{j,i}$ and so that
for each~$j$ the family~$\{W_{j,i}\}_i$ covers~$X$. Fix
definable continuous functions~$u_{j,i}\colon X\to\I$ such that
$u_{j,i}\rest{W_{j,i}}\equiv 1$ and $u_{j,i}\rest{X\setminus
V_{j,i}}\equiv 0$, which is possible by definable partition of
unity~\cite[Chapter~\n{6} Lemma~\n{3}.\n{7}]{VanDenDries}.
Finally define
\begin{align*}
	\sigma_{j,i} \colon X &\to \I \\
	x &\mapsto \max\left(0,u_{j,1}(x),\dotsc,u_{j,i}(x)\right) \\
	\tau_{j,i} \colon X &\to \I \\
	x &\mapsto g_j(x) +
		\sigma_{j,i}(x)\left(g_{j+1}(x)-g_j(x)\right)
\end{align*}
and write~$\tau_a$ for~$\tau_{j,i}$ where~$(j,i)$ is the $a$\-th pair of indices
in lexicographical order ($j$ is more important, both are increasing).

Stipulating that~$\tau_0$ denotes the constant~$0$,
we have a finite family~$\{\tau_a\}_{a\le nk}$ of continuous
definable functions with the following property: for each~$a$ let
\[
X_a=\left\{(x,t)\in X\times\I\st 0\le t\le \tau_a(x)\right\}
\]
then $X_a\subset X_{a+1}$, and for each~$a$ there is
a~$b_a\in\{1,\dotsc,n\}$ such that
\[\tilde{h}(X_{a+1}\setminus X_a)\subset U'_{b_a}\]
in fact, $b_a$ is~$i$ when $(j,i)$ is the
$a$\-th pair.
For each~$a$ let
\[
Y_a = \left\{(x,t)\in B\times\I\st \left(p(x),t\right)\in X_a\right\}
\]
The homotopy~$h$ is given on~$Y_0=B\times\{0\}$: we will extend it
inductively on the sets~$Y_\wildcard$ up to~$Y_{nk}=B\times\I$.
Fix an~$a$ and suppose to have already defined $h$ on~$Y_a$.
For $(x,t)\in Y_{a+1}\setminus Y_a$ let
\[
h(x,t) = \phi'^{-1}_{b_a}\left(\tilde{h}\left(p(x),t\right),p'_{b_a}\circ
h\left(x,\tau_a\circ p(x)\right)\right)
\]
where the formula makes sense since
\[h(x,\tau_a\circ p(x))\in\dom{p'_{b_a}}=p'^{-1}(U_{b_a})\]
which we have observing that
%$(x,t)\in Y_{a+1}\setminus Y_a$ means
%that
%\[\tau_a\circ p(x)<t\le\tau_{a+1}\circ p(x)\]
%therefore
$\{p(x)\}\times\intervaloc{\tau_a\circ p(x)}{\tau_{a+1}\circ p(x)}$ is a non\-empty subset of~$X_{a+1}\setminus X_a$, hence
$\tilde{h}(p(x),\tau_a\circ
p(x))\in\closure{U'_{b_a}}\subset{U_{b_a}}$.

Our function is clearly continuous and fibre preserving, hence a definable
bundle map. Stationarity is immediate by inspection of the formula above.
\end{proof}

\begin{corollary}\label{th-contractible-base}
Any definable fibre bundle on a definably contractible base having locally
definably compact fibre is trivial.
\end{corollary}
\begin{proof}
Immediate from Theorem~\ref{th-covering-homotopy} and
Theorem~\ref{th-equivalence-theo}.
\end{proof}

Now we give a definition of homotopy equivalence for definable
fibre bundles. Notice that two definable fibre bundles which are homotopy
equivalent have definably homotopy equivalent bundle spaces, definably
homotopy equivalent bases, and the same fibre.

\begin{definition}\label{def-homoequiv-bundles}
Two definable fibre bundles~$\bundle{B}$ and~$\bundle{B'}$ are said to
be~\new{homotopy equivalent definable fibre bundles} if there are
two definable bundle maps~$f\colon\bundle{B}\to\bundle{B'}$
and~$g\colon\bundle{B'}\to\bundle{B}$ such that $f\circ g$ and~$g\circ f$
are both homotopic to the identity on their respective domains---i.e.\ there
are homotopies \emph{of definable bundle maps} between each of them and the
identity. In this situation, $f$ and~$g$ are said to be~\new{homotopy
inverse} of each other.
\end{definition}

\begin{observation}\label{th-homotopy-inverse}
Fix an~$f\colon\square\to\triangle$ and let~$g'\colon\triangle\to\square$
be a~\new{left homotopy inverse} of~$f$, which means that
$g'\circ f$ is homotopic to the identity.
Let~$g''\colon\triangle\to\square$ be a~\new{right homotopy inverse}
of~$f$,
i.e.\ $f\circ g''\sim\id$. Then $f$ has an homotopy inverse which
is~$g'\circ f\circ g''$.
\end{observation}

\begin{theorem}\label{th-homoequiv-ind-bundle}
Let~$X$ and~$X'$ be definable sets, and let~$f\colon X\to X'$ be a
definable homotopy equivalence---i.e.\ there is a definable homotopy
inverse~$g\colon X'\to X$ of~$f$.
Consider a definable fibre bundle~$\bundle{B'}=(B',X',p,F)$
having~$X'$ as its base space and locally definably compact fibre~$F$.
Then~$\bundle{B}\eqdef f^{-1}(\bundle{B'})$ is homotopy
equivalent to~$\bundle{B'}$.
\end{theorem}
\begin{proof}
Take a definable homotopy~$\onthebase{h}\colon X'\times\I\to X'$ with $\onthebase{h}_0=\id$
and~$\onthebase{h}_1=f\circ g$. By Theorem~\ref{th-covering-homotopy} we have an
homotopy of definable bundle
maps~$h\colon\bundle{B'}\times\I\to\bundle{B'}$ such that $h_0=\id$ and
the homotopy induced by~$h$ on the base spaces coincides
with~$\onthebase{h}$---what justifies our abuse of the
notation~$\onthebase{h}$. By Theorem~\ref{th-equivalence-theo} we know that
$\bundle{B'}$ is isomorphic to~$\onthebase{h}_1^{-1}(\bundle{B'})$ which,
in turn, is isomorphic
to~$g^{-1}(\bundle{B})$ by Observation~\ref{th-composition-ind},
and we have as well a definable bundle
isomorphism~$\psi\colon\bundle{B'}\to g^{-1}(\bundle{B})$ such
that~$\indmap{f}\circ\indmap{g}\circ\psi=h_1$. Hence $\indmap{f}$ has a
right homotopy inverse, which is~$\indmap{g}\circ\psi$.

By the very same argument $\indmap{g}$ -- hence $\indmap{g}\circ\psi$ --
has a right homotopy inverse. So $\indmap{g}\circ\psi$, having both right
and left -- which is~$\indmap{f}$ -- homotopy inverses, is a definable
homotopy equivalence between~$\bundle{B'}$ and~$\bundle{B}$.
\end{proof}

\section{Bundles on the torus}\label{sect-bundles-on-tori}

\noindent In this section we will study the particular case of definable fibre
bundles whose base is a~\new{definable torus}---i.e.\ $\T^n$ for some~$n$,
where $\T$ denotes the definable group~$\intervalco{0}{1}$ with the sum modulo~$1$ and
the Pillay's topology\footnotemark. 
\footnotetext{To elude the group topology, for our purpose, we could as well have
defined~$\T$ as
$\SO_2$---i.e.\ the set~$\{(x,y)\in\Mst^2\st
x^2+y^2=1\}$ with the usual operation \emph{and the subset
topology}. Notice, however, that in general $\SO_2$
and~$\intervalco{0}{1}$ with the sum modulo~$1$ are not definably isomorphic.}
Let~$\pidef_m(\T^n)$ denote the $m$\-th o\-minimal homotopy
group of~$\T^n$, as defined in~\cite[Section~\n{4}]{BaOt09}.
By Corollary~\n{3}.\n{3} in~\cite{BeMaOt08} we have
that the higher -- i.e.\ with $m>1$ -- o\-minimal homotopy groups of~$\T^n$ are trivial.
%
%This and the results of section~\ref{sect-on-bundles} will enable us to
%give a necessary and sufficient condition for the triviality of~\emph{principal} definable
%fibre bundles on a torus.
%
The reader interested in a thorough discussion
of the o\-minimal homotopy groups is referred to~\cite{BaOt09},
however, all we need in this section is the observation below, which can be proven
straightforwardly from the definition.

\begin{observation}\label{th-pi-n-trivial}
Let~$X$ be a definably connected definable set such that~$\pidef_n(X)$ is
trivial---e.g.\ $\T^m$ for any~$m$ and for~$n>1$.
Let~$f,g\colon\I^n\to X$ be two definable maps, and suppose
that~$f\rest{\boundary{\I^n}}=g\rest{\boundary{\I^n}}$. Then $f$ and~$g$
are homotopic via an homotopy~$h\colon\I^n\times\I\to X$ such that
$h_t\rest{\boundary{\I^n}}=f\rest{\boundary{\I^n}}$~for
any~$t\in\I$.
\end{observation}

\begin{definition}\label{def-property-p}
We say that a definable set~$X$ has property~\P\ if $X$ is definably
connected and, for all~$n\in\N$,
any definable map~$f\colon\boundary{\I^n}\to X$ such that
{\def\b#1#2{\mathmakebox[\widthof{$#1$}][c]{#2}}%
\begin{equation*}\tag{$\star$}\label{eq-property-p}\begin{split}
	f(\b{x_1}{0},x_2,x_3,\dotsc,x_n) &= f(\b{x_1}{1},x_2,x_3,\dotsc,x_n) \\
	f(x_1,\b{x_2}{0},x_3,\dotsc,x_n) &= f(x_1,\b{x_2}{1},x_3,\dotsc,x_n) \\
					&\eqvdots\\
	f(x_1,x_2,x_3,\dotsc,\b{x_n}{0}) &= f(x_1,x_2,x_3,\dotsc,\b{x_n}{1})
\end{split}\end{equation*}}
is definably homotopic to a constant.
\end{definition}

\begin{observation}\label{th-equivalent-of-p}
Equivalently, a definably connected definable set~$X$ has property~\P\ if any definable map
from~$\boundary{\I^n}$ to~$X$
satisfying~\eqref{eq-property-p} extends to a map from~$\I^n$ to~$X$.
\end{observation}

\begin{lemma}\label{th-p-is-invariant}
Property~\P\ is a definable homotopy invariant---i.e.\ given two definably
homotopy equivalent definable sets~$X$ and~$X'$ one has property~\P\ if and
only if the other has.
\end{lemma}
\begin{proof}
Immediate from the definition. Take~$g\colon X\to X'$ and~$g'\colon
X'\to X$ definable maps definably homotopy inverse of each other. Suppose
that $X'$ has property~\P, and let~$f\colon\boundary{\I^n}\to X$
satisfy~\eqref{eq-property-p}. By construction $f$ is homotopic to~$g'\circ
g\circ f$, but $g\circ f\colon\boundary{\I^n}\to X'$
satisfies~\eqref{eq-property-p}, hence $g\circ f\sim\textsl{constant}$ by
hypothesis. As a consequence $f\sim g'\circ g\circ f\sim
g'(\textsl{constant}) = \textsl{constant}$.
\end{proof}

\begin{lemma}\label{th-p-implies-section}
Let~$\bundle{B}=(B,\T^n,p,F)$ be definable fibre bundle having
the $n$\-dimensional torus~$\T^n$ as the base space and a locally definably
compact definably connected definable set~$F$ as fibre. Suppose that the
bundle space~$B$ of~$\bundle{B}$ has property~\P. Then~$\bundle{B}$ admits
a definable cross section.
\end{lemma}
\begin{proof}
Let~$\tau_n\colon\I^n\to\T^n$ be the canonical map
\[
\tau_n\colon(x_1,\dotsc,x_n) \mapsto (x_1\,\mathrm{mod}\,1,\dotsc,x_n\,\mathrm{mod}\,1)
\]
By induction on the dimension~$n$, suffices to show that any definable
\emph{partial} cross section~$s$ defined on~$\tau_n(\boundary{\I^n})$ extends
to a \emph{global} cross section.

\noindent\textbf{Case $\mathbf{n=1}$.} Follows immediately from the
definable connectedness of the fibre, using the fact that
$\tau_1^{-1}(\bundle{B})$ is isomorphic to a product bundle by
Corollary~\ref{th-contractible-base}. Observe, however, that the same
argument doesn't work for~$n>1$, since we would need the fibre to be
$(n-1)$\-connected.

\noindent\textbf{Case $\mathbf{n>1}$.}
Since $B$ has property~\P, by Observation~\ref{th-equivalent-of-p}, the
map~$s\circ\tau_n$ extends to a map~$f\colon\I^n\to B$.
By Observation~\ref{th-pi-n-trivial}, $p\circ f$ is definably homotopic
to~$\tau_n$.
Since
$f\rest{\boundary{\I^n}}=s\circ\tau_n\rest{\boundary{\I^n}}$ there is a well defined definable
map~$f'\colon\T^n\to B$ such that~$f=f'\circ\tau_n$, moreover $p\circ f'$
is definably homotopic to the identity on~$\T^n$.
Hence, by
Theorem~\ref{th-covering-homotopy}, the identity on~$\bundle{B}$
is homotopic to a definable bundle map~$g\colon\bundle{B}\to\bundle{B}$
which induces~$p\circ f'$ on the
base space~$\T^n$ of~$\bundle{B}$. As a consequence, by
Theorem~\ref{th-equivalence-theo}, $\bundle{B'}\eqdef(p\circ
f')^{-1}(\bundle{B})$ is isomorphic to~$\bundle{B}$.
On the other hand, $\bundle{B'}$, which has base space
\[
B'\eqdef\left\{(x,y)\in B\times\T^n\st p(x)=p\circ f'(y)\right\}
\]
admits a definable cross section, which is~$s'\colon y\mapsto(f'(y),y)$.
It is thus proven that~$\bundle{B}$ has a definable cross section; some extra care must be
taken in order to ensure that this cross section coincides with~$s$
on~$\tau_n(\boundary{\I^n})$.

Using Observation~\ref{th-pi-n-trivial}, we get a definable homotopy
between~$p\circ f'$ and the identity which is stationary
on~$\tau_n(\boundary{\I^n})$. Hence the bundle map~$g$ given by
Theorem~\ref{th-covering-homotopy} restricts to the identity
on~$p^{-1}\circ\tau_n(\boundary{\I^n})$. Now, the bundle isomorphism
from~$\bundle{B}$ to~$\onthebase{g}^{-1}(\bundle{B})$ given by
the proof of Theorem~\ref{th-equivalence-theo} identifies $x\in\bundle{B}$
with~$(g(x),p(x))$,
which is in the graph of~$s'$ if and only
if~$f'\circ p(x)=g(x)$. However, on
$p^{-1}\circ\tau_n(\boundary{\I^n})$, we know that $g$ is the identity, so $f'\circ
p(x)=g(x)$ is equivalent to~$f'\circ p(x)=x$, which holds if and only if
$x$ is in the graph of~$s$.
\end{proof}

\section{Definably compact groups}\label{sect-compact-groups}

\noindent As proven in~\cite{Pi88} every definable
group has a unique definable manifold structure that makes it a
topological group, we assume definable groups to have the topology
induced by this structure. Notice that, using~\cite[Chapter~\n{10}
Theorem~\n{1}.\n{8}]{VanDenDries}, every definable
group is definably homeomorphic to some definable set. In this section we
will prove our main theorem -- which goes under the
designation of~\ref{th-groups-are-products} -- and show that we can not
extend it to a definable analogous of the Borel\-Scheerer\-Hofmann
splitting theorem. A few facts are needed.

\begin{fact}[{\cite[Lemma~\n{2}.\n{2}]{BeMaOt08}}]\label{th-quotient-bundle}
Let~$H$ be a definable subgroup of a definable group~$G$.
Then we have a definable fibre bundle
$\bundle{B}_{\quotient{G}{H}}\eqdef(G,\quotient{G}{H},p,H)$, where
we put on~$\quotient{G}{H}$
the quotient topology and $p\colon G\to \quotient{G}{H}$ is the projection. 
\end{fact}

\begin{fact}[{\cite[Corollary~\n{6}.\n{4}]{HruPePi08B}}]\label{th-derived-subgroup}
Let~$G$ be a definably compact definably connected group, then the derived
subgroup~$\derived{G}$ of~$G$ is a definably connected semisimple definable
group. Let~$\cccenter{G}$ denote the definably connected
component of the identity of the center of~$G$. Then
$G=\derived{G}\cdot\cccenter{G}$ and
$\Gamma_G\eqdef\derived{G}\cap\cccenter{G}$ is finite.
\end{fact}

\begin{fact}[{\cite[Theorem~\n{3}.\n{4}]{BeMaOt08}}]\label{th-homoequiv-tori}
Let~$G$ be a definably compact definably connected Abelian group
of dimension~$d$. Then~$G$ is definably homotopy equivalent to~$\T^d$.
\end{fact}

Clearly, from Fact~\ref{th-derived-subgroup}, we have that the
quotient~$\quotient{G}{\derived{G}}$
is definably isomorphic to~$\quotient{\cccenter{G}}{\Gamma_G}$. We will
call~$\bundle{B}_G\eqdef(G,\quotient{\cccenter{G}}{\Gamma_G},p,\derived{G})$
the bundle obtained from~$\bundle{B}_{\quotient{G}{\derived{G}}}$ through this
(group) isomorphism.

\begin{observation}\label{th-section-invariant}
Let~$\bundle{B'}=(B',X',p',F)$ be a definable fibre bundle, and
let~$f$ be a definable map from a definable set~$X$ to the
base~$X'$ of~$\bundle{B'}$. Consider the
bundle~$\bundle{B}\eqdef f^{-1}(\bundle{B'})$. If~$\bundle{B'}$
admits a definable cross section~$s'\colon X'\to B'$, then $\bundle{B}$ has a
definable cross section too. In fact, by definition,
\[
	B=\left\{(x,y)\in B'\times X\st p'(x)=f(y)\right\}
\]
and
\begin{align*}
	s\colon X &\to B \\
	y &\mapsto \left(s'\circ f(y),y\right)
\end{align*}
is a cross section of~$\bundle{B}$. In particular, using
Theorem~\ref{th-equivalence-theo}, we have that given two homotopy
equivalent definable fibre bundles having locally definably compact fibre,
one has a definable cross section if and only if the other has.
\end{observation}

The following lemma justifies our introduction of property~\P\ in
Section~\ref{sect-bundles-on-tori}.

\begin{lemma}\label{th-groups-have-p}
Let~$G$ be a definably connected definable group. Then~$G$ has
property~\P.
\end{lemma}
\begin{proof}
Let~$f\colon\boundary{\I^n}\to G$ be a definable map satisfying~\eqref{eq-property-p}.
Suffices to prove that $f$ is definably homotopic to a constant function.
For
each~$i\in\{1,\dotsc,n\}$ define the
projection~$p_i\colon\I^n\to\boundary{\I^n}$ by
\[
\left(p_i(x_1,\dotsc,x_n)\right)_j = \begin{cases}
	x_j	&\text{if $i\neq j$}\\
	0	&\text{if $i=j$}
\end{cases}
\]
Now let~$f_0=f$ and, for each~$i\in\{1,\dotsc,n\}$, define
\begin{align*}
f_i\colon\boundary{\I^n} &\to G \\
x &\mapsto f_{i-1}(x)\cdot\left(f_{i-1}\circ p_i(x)\right)^{-1}
\end{align*}
Each of these maps satisfies~\eqref{eq-property-p}---since $f_0$ does
and $f_i$ is defined inductively from~$f_{i-1}$ by means of operations
that preserve~\eqref{eq-property-p}. Moreover, for each~$i$, we have
$f_i\rest{S_i}\equiv\identity{G}$ where
\[
S_i\eqdef\left\{(x_1,\dotsc,x_n)\in\boundary{\I^n}\st\exists\,j\le
i\;x_j\in\{0,1\}\right\}
\]
what is easy to prove by induction. Since~$S_n=\boundary{\I^n}$, we have
that~$f_n\equiv\identity{G}$, hence suffices to prove that~$f_{i-1}\sim
f_i$ for each~$i$. We claim that~$f_{i-1}\circ p_i$ is definably homotopic
to a constant, which is enough---in fact, by the definable connectedness
of~$G$, we may assume that constant to be~$\identity{G}$, hence
\[f_{i-1}=f_{i-1}\cdot\identity{G}^{-1}\sim f_{i-1}\cdot(f_{i-1}\circ
p_i)^{-1}=f_i\]
Here is the required homotopy
\begin{align*}
\boundary{\I^n}\times\I &\to G \\
(x,t) &\mapsto f_{i-1}\circ p_i(tx)
\qedhere
\end{align*}
\end{proof}

\begin{lemma}\label{th-groups-have-sections}
Let~$G$ be a definably compact definably connect group. Then the
definable fibre bundle~$\bundle{B}_G$ admits a definable cross section.
\end{lemma}
\begin{proof}
First of all, observe that the fibre of~$\bundle{B}_G$ is definably
connected and (locally) definably compact.
By Fact~\ref{th-homoequiv-tori} we know
that~$\quotient{\cccenter{G}}{\Gamma_G}$ is definably
homotopy equivalent to~$\T^d$, for some~$d$.
Let~$f\colon\T^d\to\quotient{\cccenter{G}}{\Gamma_G}$ be a definable homotopy equivalence, and
consider the definable fibre bundle~$\bundle{B'}\eqdef
f^{-1}(\bundle{B}_G)$. By Theorem~\ref{th-homoequiv-ind-bundle} we have
that~$\bundle{B'}$ is homotopy equivalent to~$\bundle{B}_G$, hence the
respective bundle spaces are definably homotopy equivalent. Since, by
Lemma~\ref{th-groups-have-p}, the bundle space of~$\bundle{B}_G$ -- which
is~$G$ -- has property~\P, the same
holds for the bundle space of~$\bundle{B'}$. As a consequence, by
Lemma~\ref{th-p-implies-section}, $\bundle{B'}$ has a definable cross
section. The statement is thus proven using
Observation~\ref{th-section-invariant}. 
\end{proof}

\begin{theorem}\label{th-groups-are-products}
Let~$G$ be a definably compact definably connected group. Then
$G$ is definably homeomorphic to the Cartesian
product~$\derived{G}\times\quotient{\cccenter{G}}{\Gamma_G}$.
\end{theorem}
\begin{proof}
Immediate by Lemma~\ref{th-groups-have-sections}: in fact the homeomorphism
is
\begin{align*}
	\derived{G}\times\quotient{\cccenter{G}}{\Gamma_G} &\to G \\
	(x,y) &\mapsto x\cdot s(y)
\end{align*}
where $\cdot$~denotes the group operation in~$G$ and $s$~is a definable
cross section of~$\bundle{B}_G$.
\end{proof}

%The same argument proves the following generalization
%of lemma~\ref{th-principal-bundles-on-tori}.
%
%\begin{theorem}\label{th-principal-bundles-on-abelian}
%Any principal definable fibre bundle~$\bundle{B}=(B,A,p,G)$
%having a definably compact definably connected Abelian group~$A$
%as the base space is trivial if and only if its bundle space~$B$ has
%property~\P\ and $\pidef_1(A)$ acts trivially on the
%connected components of~$G$.
%\end{theorem}
%\begin{proof}[Proof (sketch).]
%Suffices to observe that the condition on the triviality of the action
%of~$\pidef_1(A)$ is invariant under homotopy
%equivalence~\emph{of definable fibre bundles}. Then follow the same
%argument in the proof of lemma~\ref{th-groups-have-sections} and
%apply lemma~\ref{th-principal-bundles-on-tori}.
%\end{proof}

The following example shows that a definably compact group may not be
(definably isomorphic to) a definable semidirect product of its derived
subgroup with some definable group.

\begin{example}\label{ex-no-semidirect-product}
We will construct a definably compact group~$G$ definable in the
pure field~$(\R,+,\cdot)$ of the real numbers such that the derived subgroup
of~$G$ has no definable semidirect complement in~$G$. From now on,
definable means definable in~$(\R,+,\cdot)$.

Let $\SO_2$ denote the unit circle~$\{a+b\i\st a^2+b^2 = 1\}$ in the
complex plane~$\R+\R\i$, equipped the group structure defined by the
complex product. Let $\T$
denote~$\intervalco{0}{1}$ with the sum modulo~$1$.
First of all, we observe that no non\-trivial homomorphism~$\T\to\SO_2$ can
be definable. This follows from the classification of one-dimensional Nash
groups in~\cite{MaSta92}, however we include here a simple argument
by contradiction. Let $\phi\colon\T\to\SO_2$ be a definable non\-trivial homomorphism.
By definable choice, we can assume without loss of generality that $\phi$
is $0$\-definable. Also, for any interval~$\intervaloo{x_0-\epsilon}{x_0+\epsilon}$,
the restriction $\phi\rest{\intervaloo{x_0-\epsilon}{x_0+\epsilon}}$
uniquely determines~$\phi$: suffices to use $\phi(Nx \mod 1) =
\left(\phi(x)\right)^N$. By o\-minimality there is an
interval~$\intervaloo{x_0-\epsilon}{x_0+\epsilon}$ over which $\phi$ is
differentiable, follows easily that $\phi$ is everywhere differentiable
(moreover $\lim_{x\to 1}\phi(x) = \phi(0) = 1$).
Now, by direct computation, $\phi'(x) = \phi(x)\phi'(0)$. Hence $\phi(x) =
e^{x\phi'(0)}$. Since $\lim_{x\to 1}\phi(x) = e^{\phi'(0)}$, follows that
$\phi'(0)$ is a multiple of~$2\pi$, non null because $\phi$ is not
trivial. This is a contradiction since $\pi$ is not algebraic.

Now let $\SU_2$ denote the
group
\[
	\SU_2 = \{a+b\i+c\j+d\k \st a^2+b^2+c^2+d^2=1\}
\]
with the usual quaternion multiplication.
Consider~$G\eqdef\quotient{(\T\times\SU_2)}{\Gamma}$ where $\Gamma$ is the
normal subgroup~$\{(0,1),(\nicefrac{1}{2},-1)\}$.
We will show that the derived subgroup of~$G$ has no definable semidirect complement
in~$G$.
It is easy to check that
$\T\times\{1,-1\}\cong\T$ is the center of~$G$, and
$\{0,\nicefrac{1}{2}\}\times\SU_2\cong\SU_2$ is the derived subgroup
of~$G$, hence $(0,-1)\cdot\Gamma$ is an element
of~$\cccenter{G}\cap\derived{G}$, which therefore is non\-trivial.
Now assume, by contradiction, that $\derived{G}$ has a definable
semidirect complement~$H<G$. Let~$\sigma\colon
G\to\derived{G}$ and~$\tau\colon G\to H$ be defined by the
equation~$x=\sigma(x)\cdot\tau(x)$ for each~$x\in G$ (these functions are
definable by definition of semidirect complement, because the decomposition is unique).
Observe that
$\sigma\rest{\gcenter{G}}$ is a group homomorphism, moreover it is the
identity on~$\gcenter{G}\cap\derived{G}$.
Since $\cccenter{G}\cap\derived{G}$ is non\-trivial,
$\sigma\rest{\cccenter{G}}$ must be a non\-trivial definable group
homomorphism from~$\cccenter{G}$ to~$\derived{G}$. 
Being a
definably connected Abelian definable subgroup of~$\SU_2$ of dimension~$1$,
the image of~$\sigma\rest{\cccenter{G}}$
must be a $0$\-Sylow in the sense of~\cite{Strze94}.
Follows that it must be conjugated to the
subgroup~$\SO_2\eqdef\{a+b\i\}<\SU_2$. However we observed that
no non\-trivial homomorphism $\T\to\SO_2$ can be definable: a
contradiction.
\end{example}

\vfill\pagebreak
\providecommand{\bysame}{\leavevmode\hbox to3em{\hrulefill}\thinspace}
\providecommand{\MR}{\relax\ifhmode\unskip\space\fi MR }
\providecommand{\MRhref}[2]{%
  \href{http://www.ams.org/mathscinet-getitem?mr=#1}{#2}
}
\providecommand{\href}[2]{#2}

\end{document}